\documentclass[runningheads]{llncs}

\usepackage{amsmath,amssymb}
\usepackage[bookmarks,
            raiselinks,
            pageanchor,
            hyperindex,
            colorlinks,
            citecolor=black,
            linkcolor=black,
            urlcolor=black,
            filecolor=black,
            menucolor=black]{hyperref}

\usepackage{multirow}

\newcommand{\mycode}[1]{\textnormal{\texttt{#1}}}

\newcommand{\fenics}{\mycode{FEniCS}}
\newcommand{\ngsolve}{\mycode{NGSolve}}
\newcommand{\deal}{\mycode{deal.II}}
\newcommand{\dunegdt}{\mycode{dune-gdt}}
\newcommand{\pymor}{\mycode{pyMOR}}
\newcommand{\R}{\mathbb{R}}
\newcommand{\N}{\mathbb{N}}
\newcommand{\ns}{\mathcal{H}_k(\mathcal{P})}
\newcommand{\dx}{\,\text{d}\hspace{-0.25pt}x}
\newcommand{\ds}{\,\text{d}\hspace{-0.25pt}s}
\DeclareMathOperator*{\argmin}{arg\min}
\DeclareMathOperator*{\argmax}{arg\max}
\DeclareMathOperator*{\Sp}{span}

\usepackage[frozencache]{minted}
\usepackage{etoolbox}
\AtBeginEnvironment{minted}{%
  \setlength{\parskip}{-3pt}
}

\begin{document}

\title{%
  A full order, reduced order and machine learning model pipeline for efficient prediction of reactive flows%
  \thanks{%
    Funded by BMBF under contracts 05M20PMA, 05M20VSA, 05M20AMD.
    Funded by the Deutsche Forschungsgemeinschaft (DFG, German Research Foundation) under Germany's Excellence Strategy EXC 2044–390685587, Mathematics Münster: Dynamics - Geometry - Structure.
    Funded by the Deutsche Forschungsgemeinschaft (DFG, German Research Foundation) under Germany's Excellence Strategy EXC 2075–390740016.
    We acknowledge the support by the Stuttgart Center for Simulation Science (SimTech).
      }
}
\titlerunning{%
  FOM, ROM and ML pipeline for efficient prediction of reactive flows%
}

\author{%
  Pavel Gavrilenko\inst{3} \and
  Bernard Haasdonk\inst{2} \and
  Oleg Iliev\inst{3} \and
  Mario Ohlberger\inst{1} \and \\
  Felix Schindler\inst{1} \and
  Pavel Toktaliev\inst{3} \and
  Tizian Wenzel\inst{2} \and
  Maha Youssef\inst{2}
}
\authorrunning{P.~Gavrilenko et al.}

\institute{%
  Mathematics Münster, Westfälische Wilhelms-Universität Münster, Einsteinstr. 62, 48149 Münster
  \email{\{mario.ohlberger,felix.schindler\}@uni-muenster.de}
\and
  Institute of Applied Analysis and Numerical Simulation, Pfaffenwaldring 57, 70569 Stuttgart
  \email{\{haasdonk,tizian.wenzel,maha.youssef-ismail\}@mathematik.uni-stuttgart.de}
\and
  Fraunhofer-Institut für Techno und Wirtschaftsmathematik, Fraunhofer-Platz 1, 67663 Kaiserslautern
  \email{\{oleg.iliev,pavel.gavrilenko,pavel.toktaliev\}@itwm.fraunhofer.de}
}

\maketitle

\begin{abstract}
We present an integrated approach for the use of simulated data from full order discretization as well as projection-based Reduced Basis reduced order models for the training of machine learning approaches, in particular Kernel Methods, in order to achieve fast, reliable predictive models
for the chemical conversion rate in reactive flows with varying transport regimes.
\end{abstract}

\section{Introduction}
\label{sec:intro}

Reactive mass transport in porous media with catalytic reactions is the basis for many industrial processes and systems, such as fuel cells, photovoltaic cells, catalytic filters for exhaust gases and catalytic burners.
The usual way of designing and testing prototypes of such devices is expensive and time-consuming.
While modeling and simulation of the processes at the pore scale of the porous media can help to optimize the design of device catalytic components, it is currently limited by the fact that such simulations lead to large amounts of data (each simulation may consist of hundreds of TB).
Moreover, processes under consideration depend on a large number of parameters.
As a consequence, the development of approaches to solve these problems with large amounts of data, as well as for the prediction of the chemical conversion rate using modern data-based methods is essential for fast, reliable predictive models.

The purpose of this paper is to
demonstrate
a computational pipeline which combines direct computational tools and different data-based methods
for a simple model problem with industry relevant aspects.
As a basic underlying model for reactive flow we consider the following scalar convection-diffusion-reaction (CDR) model problem for a concentration $c_\mu$, i.e.
\begin{align}
\label{eq:pde_0}
\partial_t c_\mu - \Delta c_\mu + \mu_2\, \nabla \cdot (u c_\mu) + \mu_1\, c_\mu &= 0,
\end{align}
posed on a unit domain with prescribed velocity $u$ of unit magnitude,
where we consider the Damköhler and Pecl\'{e}t numbers $Da, Pe \in \mathbb{R}^+$ as parameters, i.e. $\mu = (\mu_1, \mu_2)^\top = (Da, Pe)^\top$. The problem will be complemented with suitable initial and boundary data, as well as with a suitable quantity of interest that needs to be evaluated.

In Section \ref{sec:2}, we detail the problem formulation and involved concepts constituting the full approximation pipeline, including
\begin{itemize}
  \item a full order model (\textbf{FOM}),
  \item a reduced order model (\textbf{ROM}) based on FOM data and
  \item a machine learning (\textbf{ML}) based model based on FOM and ROM data.
\end{itemize}
Section 3 describes the associated software pipeline while numerical experiments in Section 4 give a proof of concept of our approach for a simple, but industrially relevant example.

\section{Approximation and data based learning }
\label{sec:2}

Given a bounded set of input parameters $\mathcal{P} \subset \mathbb{R}^p$, $p \in \mathbb{N}$ and an end time $T > 0$, we seek to efficiently and accurately approximate the evaluation of a function $f \in L^2\big(\mathcal{P}; L^2([0, T])\big)$ at a fixed set of finitely many time points for varying inputs $\mu \in \mathcal{P}$.
The function $f$, modeling a quantity of interest (QoI), is implicitly given by $f(\mu; t) := s_\mu\big(c_\mu(t)\big)$, where for any
input parameter $\mu \in \mathcal{P}$,
$s_\mu \in V'$ is a linear functional and the state $c_\mu \in L^2\big(0, T; V)$ with $\partial_t c_\mu \in L^2\big(0, T; V'\big)$
is the unique weak solution of
a parabolic partial differential equation (PDE)
\begin{align}
  \left<\partial_t c_\mu, v\right> + a_\mu(c_\mu, v) = l_\mu(v) &&\text{for all } v \in V, &&\text{ and } c_\mu(0) = c_0,
\label{eq:pde}
\end{align}
given initial data $c_0 \in V$, a Gelfand triple of suitable Hilbert spaces $V \subset H^1(\Omega) \subset L^2(\Omega) \subset V'$ associated with a spatial domain $\Omega$, a continuous linear functional $l_\mu \in V'$ and a continuous and coercive bilinear form $a_\mu: V \times V \to \R$, for each parameter $\mu \in \mathcal{P}$.

Since we consider stationary non-homogeneous Dirichlet boundary data $g_\text{D} \in H^{1/2}(\Gamma_\text{D})$ on the Dirichlet boundary $\Gamma_\text{D} \subset \partial\Omega$ and Neumann data $g_\text{N} \in L^2(\Gamma_\text{N})$ on the Neumann boundary $\Gamma_\text{N} := \partial\Omega\backslash\Gamma_\text{D}$, we select an extension of the Dirichlet data $\tilde{g}_\text{D} \in H^1(\Omega)$, such that $\tilde{g}_\text{D}|_{\Gamma_\text{D}} = g_\text{D}$ in the sense of traces, and consider the shifted solution trajectory $c_{0, \mu} := c_\mu - \tilde{g}_\text{D} \in V := H^1_{\Gamma_\text{D}}(\Omega) := \big\{v \in H^1(\Omega)\,\big|\,v|_{\Gamma_\text{D}} = 0 \text{ in the sense of traces}\big\}$.
We then obtain the weak formulation of \eqref{eq:pde_0} by noting $\partial_t \tilde{g}_\text{D} = 0$, letting $l_{0, \mu}(v) := \int_{\Gamma_\text{N}}g_\text{N}v\ds$,
$a_\mu(c, v) := \int_\Omega \big(\nabla c - \mu_2 u\, c\big)\cdot \nabla v \dx  + \int_\Omega \mu_1 c v \dx + \int_{\Gamma_\text{N}} (\mu_2 c u)\cdot n v\ds$,
and $l_\mu(v) := l_{0, \mu}(v) - a_\mu(\tilde{g}_\text{D}, v)$ in \eqref{eq:pde},
recalling $\mu = (\mu_1, \mu_2)^\top = (Da, Pe)^\top$.
Thus, \eqref{eq:pde} describes the time evolution of the shifted solution and we obtain the original solution by $c_\mu := c_{0, \mu} + \tilde{g}_\text{D}$.

\subsection{Approximation by a full order model (FOM)}
\label{sec:rb_fom}
Since we may not evaluate $f$ exactly, given any $N_T \in \mathbb{N}$, we consider a so-called full order model (FOM) approximation,
\begin{align}
  f_h: \mathcal{P} \to \mathbb{R}^{N_T}, && \mu \mapsto \big(f_h(\mu; t_1), \dots, f_h(\mu; t_{N_T})\big)^\top,\; f_h(\mu; t) := s_\mu(c_{h, \mu}(t)),
  \label{eq:f_h}
\end{align}
stemming for simplicity from an implicit Euler discretization on an equidistant temporal grid (yielding $N_T$ points $\{t_n\}_{1 \leq n \leq N_T} \subset [0, T]$),
and a conforming spatial discretization space $V_h \subset V$ of fixed dimension $N_h \in \mathbb{N}$,
yielding for each $1 < n \leq N_T$ the unique solution $c_{h, \mu}(t_n) \in V_h$ of
\begin{align}
  \tfrac{1}{\Delta t}\big(c_{h, \mu}(t_{n+1}) - c_{h, \mu}(t_n), v\big)_{L^2(\Omega)} + a_\mu(c_{h, \mu}(t_{n+1}), v) = l_\mu(v) &&\forall\, v \in V_h,
\label{eq:fom_time_step}
\end{align}
assuming $c_{h, \mu}(t_1) := c_0 \in V_h$, for simplicity.
Depending on $N_h$ and $N_T$, the computation of $c_{h, \mu}$ (and thus the evaluation of $f_h$) in the context of parameter studies, optimal design, uncertainty quantification or related applications involving $f_h$, may easily be prohibitively costly.

In the non-parametric setting, i.e.~when considering \eqref{eq:pde} for a single input tuple $\mu = (Da, Pe)^\top$, adaptive Finite Element methods
as well as adaptive time stepping schemes are the methods of choice, in particular when considering long-time integration.
While these are also applicable in the context of model order reduction, we restrict ourselves to fixed equidistant temporal and spatial grids.
As approximation space $V_h$ in \eqref{eq:fom_time_step} we use standard conforming piecewise linear finite elements, assuming for simplicity $\tilde{g}_\text{D} \in V_h$.

Note that other suitable choices include stabilized FEM or (upwind) Finite Volume (FV) schemes (in particular for large P\'{e}clet numbers which might induce steep spatial gradients or oscillations) and interior penalty discontinuous Galerkin (DG) schemes (in particular higher order variants to benefit from smooth spatial parts of the state).
\vspace{-0.5em}

\subsection{Reduced Basis reduced order model (ROM)}\vspace{-0.5em}
\label{sec:rb_rom}

Employing machine-learning (ML) techniques such as artificial neural networks or kernel methods,
one could directly utilize $f_h$ to learn an approximation
$f_\text{ml}: \mathcal{P} \to \mathbb{R}^{N_T}$ to efficiently provide cheap approximations of $f_h$.
In particular, ML-based models usually do not require the computation of a
(reduced) state and may even compute all $N_T$ values at once without time-integration - or even
provide predictions for continuous times $t$ instead of for $N_T$ discretized time values.
However, to ensure good approximation properties of such a ML-based model,
the computation of the training
and validation data involving $f_h$ may still be too demanding.
So, a method for rapidly generating additional training data for data augmentation is required.

As an intermediate step, we thus employ a structure preserving Reduced Basis (RB)
reduced order model (ROM), which we obtain by Galerkin-projection of \eqref{eq:fom_time_step} onto a carefully crafted low-dimensional RB subspace $V_\text{rb} \subset V_h$
(see \cite{MR3672144}).
While the construction of the RB space also involves solving \eqref{eq:fom_time_step}, it then allows the computation of a reduced state $c_{\text{rb}, \mu}(t_n) \in V_\text{rb}$, given for $1 \leq n \leq N_T$ as the unique solution of
\begin{align}
  \tfrac{1}{\Delta t} \big(c_{\text{rb}, \mu}(t_{n+1}) - c_{\text{rb}, \mu}(t_n), v\big)_{L^2(\Omega)} + a_\mu(c_{\text{rb}, \mu}(t_{n+1}), v) = l_\mu(v) &&\forall\, v \in V_\text{rb},
\label{eq:rb_rom_time_step}
\end{align}
and thus the evaluation of an RB QoI $f_\text{rb}: \mathcal{P} \to \mathbb{R}^{N_T}$,
\begin{align}
   \mu \mapsto \big(f_\text{rb}(\mu; t_1), \dots, f_\text{rb}(\mu; t_{N_T})\big)^\top &&\text{with}&& f_\text{rb}(\mu; t) := s_\mu(c_{\text{rb}, \mu}(t)),
   \label{eq:f_rb}
\end{align}
with a computational complexity only depending on $N_\text{rb} := \dim V_\text{rb}$, not $N_h$, owing to a pre-computation of all $V_h$-dependent quantities arising in \eqref{eq:rb_rom_time_step} under mild assumptions on the parametrization of $a_\mu$ and $l_\mu$.

Since the RB-ROM \eqref{eq:rb_rom_time_step} simply arises as the Galerkin projection of \eqref{eq:pde} onto $V_\text{rb}$, it is fully defined once we specify the RB space, the
construction of which is a delicate matter: it should be as low-dimensional as possible, to ensure a good online-performance of the resulting RB model;
it should be as rich as possible, to ensure good approximation properties of the resulting RB model;
however, at the same time, its construction greatly impacts the overall performance of the scheme and should be as cheap as possible.
For simplicity we employ the method of snapshots
: we collect a series of state trajectories $\{c_{h, \mu}\}_{\mu \in \mathcal{P}_\text{POD}}$ for a set of a priori specified training parameters $\mathcal{P}_\text{POD} \subset \mathcal{P}$ and simply obtain $V_\text{rb}$ by proper orthogonal decomposition (POD) of the resulting snapshot Gramian.

Note that the approximation quality of the resulting model can only be assessed in a post-processing step involving the computational complexity of the FOM again.
An alternative adaptive construction of $V_\text{rb}$ can be achieved by means of an iterative POD-greedy algorithm steered by a posteriori error control, which at the same time allows an efficient quantification of the induced model reduction error independent of the FOM complexity.
However, the computational complexity of the required offline-computations may be significant (in particular for problems with long time-integration).

\medskip

In particular for spatially highly resolved FOMs with $N_h \gg 1$, evaluating $f_\text{rb}$ might be orders of magnitude faster than $f_h$.
However, the solution of \eqref{eq:rb_rom_time_step} still requires time integration which is why for temporally highly resolved FOMs with $N_T \gg 1$,
we employ advanced machine learning ROMs using greedy kernel methods to
directly learn the mapping $f_\text{ml}: \mathcal{P} \rightarrow \R^{N_T}$, thus skipping the time integration (detailed in the following section).

\subsection{Approximation by machine learning: Kernel methods}\vspace{-0.2em}
\label{sec:ml_rom}

We consider $\mathcal{P} \subset \R^2$ and for our problem complexity shallow instead of deep learning
architectures are sufficient.
We apply kernel methods which in our case rely on strictly positive definite kernels, which are in the scalar case symmetric functions $k: \mathcal{P} \times \mathcal{P} \rightarrow \R$, such that the so called kernel matrix $\left( A_n \right)_{i,j} = k(\mu_i, \mu_j)$, $i,j = 1, .., n$ is positive definite for any set $\{ \mu_1, .., \mu_n \} \subset \mathcal{P}, n \in \N$ of pairwise distinct inputs (see \cite{Wendland2005}). The well known Gaussian kernel $k(x,y) = \exp(-\Vert x - y \Vert_2^2)$ is an example of a strictly positive definite kernel.

Associated to a strictly positive definite kernel on a domain $\mathcal{P}$, there is a space of continuous function, the so called Reproducing Kernel Hilbert Space (RKHS), $\ns$. For given input data $\{\mu_1, .., \mu_N \} \subset \mathcal{P}$ and corresponding target data $\{y_1, .., y_N \}$ $\subset \R^{d_\text{out}}$, learning with kernels refers to a minimization of a loss-functional over the RKHS. Using a mean-squared error loss with a standard norm regularization, a kernel representer theorem states
that a solution
{%
\setlength{\abovedisplayskip}{3pt}
\setlength{\belowdisplayskip}{3pt}
\begin{align}
\label{eq:loss}
f_\text{ml} := \argmin_{f \in \ns} \mathcal{L}(f), \quad \mathcal{L}(f) = \frac{1}{N} \sum_{i=1}^N \Vert y_i - f(\mu_i) \Vert_{\R^{d_\text{out}}}^2  + \lambda \cdot \Vert f \Vert_{\ns}^2
\end{align}
}%
can be found in the finite dimensional subspace spanned by the data, i.e.\ there exists $f_\text{ml} \in \Sp \{ k(\cdot, \mu_i), i = 1, .., N \}$. In order to learn a sparse and fast surrogate $f_\text{ml}$, one strives not to use all the training data, but only a meaningful subset. While a global optimization is combinatorial infeasible, we use greedy methods in order to select a subset, as implemented in the vectorial kernel orthogonal greedy algorithm (VKOGA) \cite{Santin2021}. Algorithms of this type start with an empty set $X_0 := \{ \}$, and incrementally select the next point $X_{n+1} := X_n \cup \{ x_{n+1} \}$, $x_{n+1} \in \{ \mu_1, .., \mu_N \}$, for instance according to the $f$-greedy selection criterion,
$x_{n+1} := \argmax_{\mu_i, i=1, .., N} \Vert (y_i-s_n)(\mu_i) \Vert_{\R^{d_\text{out}}},$
whereby $s_n$ is the kernel approximant based on the input data $X_n$ and corresponding target data. These greedy kernel approximation algorithms have been studied thoroughly in terms of approximation speed, stability and generalization, for instance the unregularized case, $\lambda = 0$ within \eqref{eq:loss}, in \cite{wenzel_fgreedy}.

\section{Software environment}
\label{sec:software}

We aim for a flexible, user-friendly and performant software environment to meet the requirements arising from Section \ref{sec:2}, with the Python-based and freely available model reduction library \pymor\footnote{%
  Available at \url{https://pymor.org}, including references to other software libraries.
} \cite{MR3565558} at its core:
\\\textbf{Interchangable FOM:} As argued in Section \ref{sec:rb_fom}, FV schemes or conforming or DG FEM schemes might be required for the problem at hand, the choice of which is not always clear a priori.
Thus, while \pymor{}'s built-in \texttt{numpy/scipy}-based discretization might be desirable for quick prototyping, more advanced multi-purpose discretization libraries such as \deal, \dunegdt\footnote{%
  Available at \url{https://github.com/dune-community/dune-gdt}.
}, \fenics{} or \ngsolve{} (all freely available) are often required for more advanced problems.
Some applications, for instance for multi species reactive porous media transport with non-linear source terms on computer tomography-based geometries, require more specialized libraries such as Fraunhofers in house library \texttt{PoreChem} \cite{Greiner2019}.

Convenient \pymor{} wrappers are available for all of the above libraries with a well-defined API, to allow a unified handling of the resulting FOM, regardless of its origin.
For instance, given the model problem described in Section \ref{sec:2} with scalar Dammköhler- and P\'{e}clet-number as input parameters, we can call
\begin{minted}[texcomments,fontsize=\footnotesize,numbersep=8pt]{python}
c_h = fom.solve({'Da': 1, 'Pe': 1e-3})   # compute $c_{h, \mu}$ from \eqref{eq:fom_time_step}
f_h = fom.output({'Da': 1, 'Pe': 1e-3})  # compute $f_h(\mu)$ from \eqref{eq:f_h}
\end{minted}
on any FOM obtained from one of the above libraries\footnote{%
  The code examples throughout this section contain actual (shortened) Python user-code, usually encountered in \pymor{}-based applications.
}.
\\\textbf{Interchangable MOR algorithms:} As argued in Section \ref{sec:rb_rom}, there exist several established algorithms for the generation of a reduced basis, most prominently (for instationary problems) the POD (or method of snapshots) and the POD-greedy algorithm, where the applicability of each algorithm in terms of accuracy and performance is often not clear a priori.
Once wrapped as \pymor{} models, all FOMs stemming from the above libraries expose access to (the application of) their operators (e.g.~those induced by $a_\mu$ and $l_\mu$, products such as the $L^2$- or $H^1$-semi product) and all vector arrays (containing state snapshots such as $c_\mu$) of these FOMs fulfill the same API.
\pymor{} thus ships a large variety of generic algorithms applicable such as a stabilized \mycode{gram\_schmidt}, an \mycode{rb\_greedy} and \mycode{rb\_adaptive\_weak\_greedy} and \mycode{pod} algorithm, as well as the distributed or incremental \mycode{hapod} algorithm (see below).
For instance, calling
\begin{minted}[texcomments,fontsize=\footnotesize,numbersep=8pt]{python}
snapshots = fom.solution_space.empty()
for mu in parameter_space.sample_randomly(10):
    snapshots.append(fom.solve(mu))
pod_basis, _ = pod(snapshots, product=fom.h1_0_product)
\end{minted}
computes a POD basis for any FOM stemming from one of the above libraries.
Similarly, the Galerkin projection of FOMs (i.e.~its operators, products and outputs) is provided generically, yielding a structurally similar ROM as in:
\begin{minted}[texcomments,fontsize=\footnotesize,numbersep=8pt]{python}
pod_reductor = InstationaryRBReductor(
    fom, RB=pod_basis, product=fom.h1_0_product)
pod_rom = pod_reductor.reduce()
\end{minted}
\textbf{Expandable and user-friendly environment:} The high-level interactive nature of the Python programming language allows quick prototyping for beginners and experts alike, while its rich and free ecosystem allows access to high-performance libraries such as the above mentioned ones or \mycode{pytorch}\footnote{\url{https://pytorch.org/}}, often used in the context of machine learning.
In addition, \pymor{}s API and generic algorithms allows for flexible high-level user code: since all models behave similarly, and all outputs and reduced data structures are \mycode{numpy}-based, a call like
\begin{minted}[texcomments,fontsize=\footnotesize,numbersep=8pt]{python}
f_rb = pod_rom.output({'Da': 1, 'Pe': 1e-3})     # $f_{rb}(\mu)$ from \eqref{eq:f_rb}
abs_linf_output_err = np.max(np.abs(f_h - f_rb)) # $\|f_h(\mu) - f_\textnormal{rb}(\mu)\|_{L^\infty}$
\end{minted}
works for any combination of FOM and generated ROM.
It is thus easily possible to prototype and evaluate new algorithms (such as advanced MOR and ML algorithms) as well as create custom applications, such as the workflow from Section \ref{sec:2}, involving the greedy kernel methods from the VKOGA library \cite{Wirtz2013}.

\section{Numerical experiments: Reactive flow}
\label{sec:experiments}

\label{sec:model_problem}

To demonstrate the approach detailed in Section \ref{sec:2}, we consider
the one dimensional CDR equation for dimensionless molar concentration variable, $c_{\mu}(x,t)$ $\in$ $\mathbb{R}^+$, i.e.~\eqref{eq:pde_0} posed on the unit domain $\Omega := (0, 1)$ with $u = 1$, $c_0 := 0$ initial values, Dirichlet boundary $\Gamma_\text{D} := \{0\}$, Dirichlet data $g_\text{D} = 1$, Neumann boundary $\Gamma_\text{N} := \{1\}$ and Neumann data $g_\text{N} = 0$.
Here we choose diffusion time as a characteristic time and $T = 3$ to ensure a near-stationary QoI, namely the break through curve $s_\mu(t) := \int_{\Gamma_\text{N}}c_\mu(t)\ds$, at the end of the simulation.
This model
is widely used as a basis in chemical engineering and industry plug-flow/perfectly stirred reactors models, and as a consequence part of real designing processes in the industry.
It is
an excellent compromise between the complexity of real industrial models of the catalytic process, and a simple mathematical formulation providing main features of transport processes.
For this choice of initial- and boundary data, an analytical solution of \eqref{eq:pde_0} is available owing to \cite{Genuchten1982}.

We consider the diffusion dominated regime, i.e.~$\mu = (Da, Pe)^\top \in \mathcal{P} := [10^{-3}, 1]^2$ and
are interested in an overall relative approximation error w.r.t.~the target QoI $f$ in $L^\infty\big(\mathcal{P}_\text{test}; L^2([0, T])\big)$ of less than one percent, measured over a finite test set $\mathcal{P}_\text{test} \subset \mathcal{P}$.
We thus require a relative error of $10^{-4}$ from all approximation components, and use the analytical solution from \cite{Genuchten1982} to calibrate the FOM to compute the reference QoI $f_h$ with a relative error less than $10^{-4}$, yielding $h = 2^{-6}$ and $\Delta t = 2^{-13}$ (thus $N_T = 24576$ time steps) as sufficient for the diffusion dominated regime.
We also use a finer spatial grid in a second experiment to additionally represent the state $c_\mu$ accurately.
As spatial product over $V$ we chose the full $H^1$-product, which is also used for orthonormalization.
For the discretization, we use the preliminary Python bindings of \dunegdt{}\footnote{Similar to \url{https://github.com/ftschindler-work/dune-gdt-python-bindings}.} to provide a \pymor{}-compatible FOM, as detailed in Section \ref{sec:software}.

For the method of snapshots, we select the four outermost points of $\mathcal{P}$ as training parameters $\mathcal{P}_\text{POD}$ and use \pymor{}'s implementation of the incremental Hierarchical approximate POD (HAPOD) from \cite{HLR2018}, \mycode{inc\_vectorarray\_hapod}, with a tolerance of $10^{-4}$.
Handling a dense snapshot Gramian of size $(4\cdot 24576)^2$, as in the classical POD, would otherwise be infeasible.
We use a VKOGA implementation\footnote{Available at \url{https://github.com/GabrieleSantin/VKOGA}.} with default parameters and train it using the four inputs $\mu \in \mathcal{P}_\text{POD}$ and already computed FOM outputs $f_h(\mu)$ used to build the POD basis, as well as 196 randomly selected inputs $\mu \in \mathcal{P}$ and corresponding RB-ROM approximations $f_\text{rb}(\mu)$.

The performance and approximation properties of the resulting approximate models\footnote{Computed on a dual socket compute server equipped with two Intel Xeon E5-2698 v3 CPUs with 16 cores running at 2.30GHz each and 256GB
of memory available.} are given in Table \ref{tab:results}.

\begin{table}[t]
  \centering%
  \caption{%
    Accuracy and runtime (in seconds) of the FOM, RB-ROM and VKOGA models (and respective approximations of $f$) from the proposed pipeline for the experiment from Section \ref{sec:experiments} for an ``output-accurate'' FOM (first row) and a ``state-and-output-accurate'' FOM (second row).
    Offline time comprises all parts to build the model (FOM: grid + FEM assembly; RB-ROM: FOM + snapshots + HAPOD + Galerkin projection onto $V_\text{rb}$; VKOGA: RB-ROM + snapshots + fitting).
    Online time denotes average time to evaluate the respective $f_*$ (FOM: four inputs; RB-ROM: ten inputs; VKOGA: 1000 inputs).
    ``rel.~err.'' denotes respective relative errors (RB-ROM: $\|f_h(\mu) - f_\text{rb}(\mu)\|_{l_2}/\|f_h(\mu)\|_{l_2}$ over $\mathcal{P}_\text{POD}$ + five random inputs; VKOGA: $\|f_\text{rb}(\mu) - f_\text{ml}(\mu)\|_{l_2}/\|f_\text{rb}(\mu)\|_{l_2}$ over 50 random inputs).
    ``p.~o.'' denotes pay-off of the full pipeline after this many queries of $f_\text{ml}$ compared to $f_h$ (respective offline + online).
  }
  \label{tab:results}
  \vspace{-3mm}
  \begin{tabular}{c|c|c||c|c|c|c||c|c|c|c|c}
    \multicolumn{3}{c||}{FOM ($f_h$)} & \multicolumn{4}{c||}{RB-ROM ($f_\text{rb}$)} & \multicolumn{4}{c|}{VKOGA model ($f_\text{ml}$)}\\\hline\hline
    $N_h$ & \,\;offline\;\, & \,online\, & \,offline\, & $N_\text{rb}$ & \,online\,  & rel.~err. & \,offline\, & points & \,online\,  & rel.~err. & p.~o.\\\hline
    65 & 3.62e-2 & 6.99e0 & 1.56e2 & 12 & 4.15e0 & 7.81e-6 & 9.76e2 & 51 & 4.12e-4 & 3.31e-6 & 140\\\hline
    65537 & 1.63e0 & 3.39e3 & 2.47e4 & 14 & 4.21e0 & 3.31e-6 & 2.55e4 & 51 & 4.15e-4 & 3.19e-6 & 8
  \end{tabular}\vspace{-1em}
\end{table}

\newpage
\section{Conclusion}

We propose a pipeline of consecutively built full order, reduced order and ma\-chi\-ne-learned models to approximate the evaluation of a QoI depending on a parabolic PDE, and numerically analyze its performance for an industrially relevant one-dimensional problem (see Tab. \ref{tab:results}).
While the similar dimensions of the RB-ROM, $N_\text{rb}$, and the number of selected VKOGA points for both runs indicate, that the low spatial resolution for the state is sufficient to approximate the QoI in this example, the proposed pipeline pays off after 140 queries, compared to only using the FOM.
The second run demonstrates an even more pronounced pay off after 8 queries for higher spatial resolutions (as in higher dimensions).
\enlargethispage{\baselineskip}

\end{document}